\documentclass[11pt]{amsart}
\usepackage{amsmath,amsfonts,latexsym,amssymb,amsthm,mathrsfs,array}

\usepackage{amsfonts}
\usepackage[usenames]{color}
\usepackage[colorlinks=true,linkcolor=webgreen,filecolor=webbrown,citecolor=webgreen]
{hyperref}
\definecolor{webgreen}{rgb}{0,.5,0}
\definecolor{webbrown}{rgb}{.6,0,0}

\newcommand{\CC}{{\mathbb C}}
\newcommand{\ZZ}{{\mathbb Z}}

\def\gg{{\mathfrak{g}}}
\def\hh{{\mathfrak{h}}}
\def\CC{{\mathbb C}}
\def\ZZ{{\mathbb Z}}

\newtheorem{dfn}{Definition}[section]
\newcommand{\bdfn}{\begin{dfn}\rm}
\newcommand{\edfn}{\end{dfn}}
\newtheorem{thm}[dfn]{Theorem}
\newcommand{\bthm}{\begin{thm}}
\newcommand{\ethm}{\end{thm}}
\newtheorem{lmma}[dfn]{Lemma}                   
\newcommand{\blmma}{\begin{lmma}}                   
\newcommand{\elmma}{\end{lmma}}                   
\newtheorem{ppsn}[dfn]{Proposition}
\newcommand{\bppsn}{\begin{ppsn}}
\newcommand{\eppsn}{\end{ppsn}}
\newtheorem{crlre}[dfn]{Corollary}
\newcommand{\bcrlre}{\begin{crlre}} 
\newcommand{\ecrlre}{\end{crlre}}
\newtheorem{rmk}[dfn]{Remark}
\newcommand{\brmk}{\begin{rmk}\rm} 
\newcommand{\ermk}{\end{rmk}}

\numberwithin{equation}{section}

\title[Integrable modules for graded Lie tori]{Integrable modules for graded Lie tori with finite-dimensional weight spaces}
\author{Souvik Pal}
\address{Souvik Pal, Harish-Chandra Research Institute (HBNI), Chhatnag Road, Jhunsi, Prayagraj(Allahabad) 211019, Uttar Pradesh, India}
\email{pal.souvik90@gmail.com, souvikpal@hri.res.in}
\date{}

\begin{document}

\maketitle
\begin{abstract}
An important problem in the representation theory of affine and toroidal Lie algebras is to classify all possible irreducible integrable modules with finite-dimensional weight spaces. Recently the irreducible integrable modules having finite-dimensional weight spaces with non-trivial central action have been classified for a more general class of Lie algebras, namely the graded Lie tori. In this paper, we classify all the irreducible integrable modules with finite-dimensional weight spaces for this graded Lie tori where the central elements act trivially. Thus we ultimately obtain all the simple objects in the category of integrable modules with finite-dimensional weight spaces for the graded Lie tori.

\end{abstract}

\section{Introduction}

One of the main ingredients that make the study of affine Kac-Moody algebras and their representations tractable is the existence of an explicit realization. The classical procedure to acquire this concrete realization in case of a twisted affine Lie algebra proceeds in two simple steps \cite{K}. In the first step, the derived Lie algebra of the affine Kac-Moody algebra modulo its centre is constructed by means of the loop algebra (in one variable) twisted by a diagram automorphism of a finite-dimensional simple Lie algebra. In the second step, the affine Lie algebra is ultimately obtained from this loop algebra by taking its universal central extension and then adding the degree derivation. The multivariable generalization of this classical object is the subject of the current paper.

In order to generalize the above two step procedure, we first need to find a suitable replacement for the twisted loop algebra. This eventually leads us to the notion of a $centreless$ $Lie$ $torus$ satisfying $fgc$ $condition$ (which means that the Lie torus is finitely generated as a module over its centroid). It was shown in \cite{ABFP} that such a Lie torus of non-zero nullity always has a $multiloop$ $realization$ and in the one variable case, this Lie algebra is exactly equivalent to the twisted version of the loop algebra mentioned above. It is worth mentioning here that centreless Lie tori satisfying fgc condition play an extremely important role in the theory of extended affine Lie algebras (EALAs for short) as they turn out to be the $centreless$ $cores$ $of$ $almost$ $every$ EALA \cite{Ne}.

We first start with a centreless Lie torus $LT$ (see Definition \ref{Tori}) and then consider its universal central extension given by $\overline{LT}= LT \oplus \mathcal{Z}(\underline{m})$ (see Subsection \ref{Extension}) which is naturally $\mathbb{Z}^n$-graded. Finally we construct the graded Lie torus $\widetilde{LT}$ by adding the space $D$ spanned by only finitely many degree derivations to $\overline{LT}$ (see \eqref{Derivations}). The centre of this Lie algebra is spanned by $n$ elements. If all these central elements act trivially on an $\widetilde{LT}$-module, we say that the $level$ of this representation is $zero$, otherwise we say that the representation has $non$-$zero$ $level$.

Unlike their structure, the representations of the graded Lie tori are not so well-understood except in some isolated cases like affine Lie algebras (both twisted and untwisted) and untwisted toroidal Lie algebras. The study of representations of $\widetilde{LT}$ is nothing but the study of $\mathbb{Z}^n$-graded representations of the Lie torus $\overline{LT}$. Perhaps the most interesting and widely studied class of representations in case of both affine and toroidal Lie algebras are the integrable ones as they can be lifted to their corresponding groups.

The main purpose of our paper is to classify all the $level$ $zero$ irreducible integrable modules for $\widetilde{LT}$ having finite-dimensional weight spaces. If $n=1$, an exhaustive collection of these irreducible modules can be deduced from the works of V. Chari, A. Pressley and S. Eswara Rao for the untwisted as well as twisted cases \cite{C,CP2,CP1,CP3,E1}. In the multivariable setup, if the $n$-tuple of automorphisms $\sigma=(\sigma_1, \cdots, \sigma_n)$ corresponding to our Lie torus $L(\mathfrak{g}, \sigma)$ (see Definition \ref{D2.1}) are all taken to be trivial, then $\widetilde{LT}$ is simply the $untwisted$ toroidal Lie algebra. In this context, the aforementioned irreducible modules have been classified independently by S. Eswara Rao \cite{E2} and Y. Yoon \cite{YY}. Subsequently another related problem was also solved for a $specific$ $twisted$ multiloop algebra in \cite{BT}. In this particular instance, the authors did not consider the universal central extension and provided a classification only for the $graded$ $version$ of the specific multiloop algebra.

This paper is an attempt to generalize the aforesaid classification problems for the first time in case of $level$ $zero$ modules. The ideas utilized in this paper are inspired from the toroidal case \cite{E2, E3, E4}. But there are many instances where we need to resort to completely different arguments which will be evident in due course. We now summarize the content of our paper.

In Section \ref{Notation}, we construct the graded Lie torus $\widetilde{LT}$ as mentioned above. This Lie algebra has a natural triangular decomposition (see (\ref{LT})) given by
\begin{align*}
\widetilde{LT} = \big (\widetilde{LT}\big )_{-} \oplus \big (\widetilde{LT}\big )_{0} \oplus \big (\widetilde{LT}\big )_{+}.
\end{align*}
In Section \ref{Representations}, we prove the existence of a highest weight vector in our level zero irreducible integrable module $V$ (see Proposition \ref{Highest}) with respect to this triangular decomposition. In Section \ref{Centre}, we consider the highest weight space and understand the action of $\big (\widetilde{LT}\big )_{0}$ on this space, thereby showing that the action of $\mathcal{Z}(\underline{m})$ is trivial on $V$ (see Theorem \ref{Trivial}). The main novelty in this classification problem lies in Section \ref{Non-graded} where we need to employ completely different techniques to the ones used in the toroidal case to pass on to an irreducible finite-dimensional module for the multiloop algebra $LT$ (see Theorem \ref{Passing}). In this case, the main hindrance to earlier approaches is that the subalgebra of $\big (\widetilde{LT}\big )_0$ excluding the degree derivations modulo $\mathcal{Z}(\underline{m})$  $is$ $not$ $necessarily$ $abelian$ (see Remark \ref{Non-abelian}). In Section \ref{Graded}, we construct the $shifted$ $loop$ $module$ (see (\ref{Loop})) and thereby infer that our original module $V$ is actually isomorphic to an irreducible component of this loop module (see Theorem \ref{Reducible}).

In Section \ref{Final}, we recollect some results from \cite{L} which classify all the finite-dimensional irreducible representations of $LT$ (see Proposition \ref{Multiloop}). Using this classification theorem, we finally conclude that $every$ $irreducible$ $integrable$ $module$ $over$ $\widetilde{LT}$ $of$ $level$ $zero$ $comes$ $from$ $a$ $finite$- $dimensional$ $irreducible$ $representation$ $of$ $a$ $finite$-$dimensional$ $semisimple$ $Lie$ $algebra$ (see  Theorem \ref{Classification}). Such finite-dimensional irreducible representations are completely parametrized by the set of all dominant integral weights.

In a fairly recent paper by S. Eswara Rao and S. Sharma, the authors classified all the irreducible integrable representations of $\widetilde{LT}$ of $non$-$zero$ $level$ \cite{ES}. Combining their main result \cite[Theorem 9.1]{ES} along with our final classification theorem, we therefore obtain a complete list of all possible irreducible integrable modules for the graded Lie torus.\\\\
$\bf{Acknowledgements.}$ I am grateful to Prof. S. Eswara Rao for suggesting the problem. I also thank him for a careful reading of the entire manuscript and pointing out some errors in an earlier draft. Finally I record my sincere gratitude to the anonymous referee for useful comments and suggestions. 

\section{Notations and Preliminaries}\label{Notation}
Throughout the paper, all the vector spaces, algebras and tensor products are over the field of complex numbers $\mathbb{C}$. We shall denote the set of integers, natural numbers and non-negative integers by $\mathbb{Z}$, $\mathbb{N}$ and $\mathbb{Z_{+}}$ respectively. For any Lie algebra $L$, the centre of $L$ will be denoted by $Z(L)$ and its universal enveloping algebra will be denoted by $U(L)$.
\subsection{Toroidal Lie Algebras}\label{Toroidal}
Let us consider a finite-dimensional simple Lie algebra $\mathfrak{g}$ equipped with a Cartan subalgebra $\mathfrak{h}$. Then $\mathfrak{g}$ is endowed with a symmetric, non-degenerate bilinear form which is invariant under every automorphism of $\mathfrak{g}$. Let us denote this form by $(\cdot|\cdot)$. Moreover we also know that the root space decomposition of $\mathfrak{g}$ automatically induces a triangular decomposition of $\mathfrak{g}$ given by $\mathfrak{g} = \mathfrak{n}^{-} \oplus \mathfrak{h} \oplus \mathfrak{n}^{+}$.\\ 
Let $A = \CC[t_1^{\pm 1}, \cdots, t_n^{\pm 1}]$  be the algebra of Laurent polynomials in $n$ variables.
Consider the (untwisted) multiloop algebra given by
\begin{equation}\label{Algebra}
L(\mathfrak{g})=\mathfrak{g} \otimes A 
\end{equation}
which is a Lie algebra under the following pointwise bilinear operation.
\begin{align*}
[x \otimes f, y \otimes g] = [x,y] \otimes fg \ \forall \ x,y \in \mathfrak{g} \ \text{and} \ f, g \in A.
\end{align*}
For any $x \in \mathfrak{g}$ and $\underline{k} \in \mathbb{Z}^n$, write $t^{\underline{k}}=t_1^{k_1} \cdots t_n^{k_n}$ and let $x(\underline{k}) = x \otimes t^{\underline{k}}$ denote a typical element of $L(\mathfrak{g})$. Now consider the module of differentials ($\Omega_A,d$) of $A$ which is the free $A$-module with basis $\{K_1, \cdots, K_n \}$ along with the differential map $d : A \longrightarrow \Omega_A$ whose image is the subspace spanned by the elements of the form $\sum_{i =1}^{n}k_i t^{\underline{k}}K_i$ for $\underline{k} \in \ZZ^{n}$. More precisely, we have       
\begin{equation}\label{Subspace}
\Omega_A = \text{span} \{t^{\underline{k}}K_i \ | \ 1 \leq i \leq n, \ \underline{k} \in \ZZ^{n} \},\ 
 dA =  \text{span} \{\sum_{i =1}^{n}{k_i t^{\underline{k}}K_i} \} \subseteq \Omega_A.
\end{equation}
If we now consider the quotient space $\mathcal{Z} = \Omega_A / d_A$, then we know that
\begin{align*}
\overline{L}(\mathfrak{g})=L(\mathfrak{g}) \oplus \mathcal{Z}
\end{align*} 
is the universal central extension of $L(\mathfrak{g})$ \cite{Ka,EMY}. By abuse of notation, we shall denote the image of $t^{\underline{k}}K_i$ in $\mathcal{Z}$ again by itself. Note that both  $L(\mathfrak{g})$ and $\overline{L}(\mathfrak{g})$ are naturally $\mathbb{Z}^n$-graded. To keep track of this grading, let us introduce the derivations $d_1, \cdots, d_n$ on $\overline{L}(\mathfrak{g})$ (whose actions are defined below) and take $D$ to be the linear span of $d_1, \cdots ,d_n$. Finally if we define 
\begin{align*}
\tau_A = L(\mathfrak{g}) \oplus \mathcal{Z} \oplus D,
\end{align*} 
then $\tau_A$ is called the (untwisted) toroidal Lie algebra (corresponding to $\mathfrak{g}$) along with the following bracket operations. 
\begin{enumerate}
\item $[x(\underline{k}), y(\underline{l})] = [x,y](\underline{k} + \underline{l}) + (x|y) \sum_{i=1}^{n} k_it^{\underline{k} + \underline{l}}K_i$;
\item $\mathcal{Z}$ is central in $\overline{L}(\mathfrak{g})$;
\item $d_i(x(\underline{k}))=[d_i, x(\underline{k})] = k_i x(\underline{k}) $;
\item $d_i(t^{\underline{m}}K_j)=[d_i, t^{\underline{m}}K_j]= m_it^{\underline{m}}K_j$;
\item $[d_i, d_j] = 0 \ \forall\ i,j = 1, \cdots ,n$.
\end{enumerate}
This Lie algebra clearly admits of a natural triangular decomposition
\begin{equation}\label{Abelian}
\tau_A = \tau_{A}^{-} \oplus \tau_{A}^{0} \oplus \tau_{A}^{+}
 \end{equation}
 where
\begin{align*}
\tau_{A}^{-} =  \mathfrak{n}_{-} \otimes A, \,\,\,\,\,\,\,\,\ 
\tau_{A}^{+} = \mathfrak{n}_{+} \otimes A, \,\,\,\,\,\,\,\,\
\tau_{A}^{0} = \mathfrak{h} \otimes A \oplus \mathcal{Z} \oplus D.
\end{align*}

\subsection{Twisted Multiloop Algebras and Lie Tori}\label{Twisted}
Fix any $n\in\mathbb{N}$ and suppose that we have $n$ commuting automorphisms $\sigma_1, \cdots, \sigma_n$ of $\mathfrak{g}$  given of finite orders $m_1, \cdots, m_n$ respectively. Put
\begin{align*} 
\sigma = (\sigma_1, \cdots, \sigma_n),\,\,\,\,\,\,\,\
\Gamma = m_1 \mathbb{Z} \oplus \cdots \oplus m_n\mathbb{Z},\,\,\,\,\,\,\,\
 \overline{\Lambda} = \mathbb{Z}^n/ \Gamma.
 \end{align*} 
Thus we have a natural map \ $\mathbb{Z}^n \longrightarrow \overline{\Lambda}$ $(\cong \mathbb{Z}/m_1\mathbb{Z} \times \cdots \times\mathbb{Z}/m_n\mathbb{Z})$ 
\begin{align*}
(k_1, \cdots, k_n) = \underline{k} \mapsto \overline{k}=(\overline{k_1}, \cdots ,\overline{k_n})
\end{align*} 
For $1 \leqslant i \leqslant n$, let $\xi_i$ denote a $m_i$-th primitive root of unity.\\
Then we obtain an eigenspace decomposition of $\mathfrak{g}$ given by
\begin{align*}
\mathfrak{g} = \bigoplus _{\overline{k} \in \overline{\Lambda}} \mathfrak{g}_{\overline{k}} \,\,\,\ \text{where} \,\,\,\ \mathfrak{g}_{\overline{k}} = \{x \in \mathfrak{g} \ | \ \sigma_i x = \xi_i^{k_i}x,\ 1 \leqslant i \leqslant n \}.
\end{align*}
It is well-known that $\mathfrak{g}_{\overline{0}}$ is a reductive Lie algebra  \cite[Proposition 4.1]{BM}, even with the possibility of being zero.
Finally let us define
\begin{align*}
L(\mathfrak{g}, \sigma) = \bigoplus_{\underline{k} \in \mathbb{Z}^n} (\mathfrak{g}_{\overline{k}} \otimes \mathbb{C}t^{\underline{k}})
\end{align*}
which is clearly a Lie subalgebra of the (untwisted) multiloop algebra $L(\mathfrak{g})$.
\begin{dfn}\label{D2.1}
The above Lie algebra $L(\mathfrak{g}, \sigma)$ is known as the (twisted) multiloop algebra associated to $\mathfrak{g}$ and $ \sigma$.
\end{dfn}
For any finite-dimensional simple Lie algebra $\gg_1$ equipped with a Cartan subalgebra $\hh_1$ of $\gg_1$, we have the root space decomposition of $\gg_1$ relative to $\hh_1$ given by
\begin{align*}
\gg_1 = \bigoplus_{\alpha \in \hh_1^*}\mathfrak{g}_{1, \alpha}
\end{align*}
$\,\,\,\,\,\,\,\,\,\,\,\,\,\,\,\,\,\,\,\,\,\,\,\,\,\,\,\,\,\,\,\  \text{where} \ \mathfrak{g}_{1, \alpha} = \{ x \in \gg_1 \ | \ [h,x]= \alpha(h)x \ \forall \ h \in \hh_1\}$.\\
Set $\Delta(\gg_1, \hh_1) = \{ \alpha \in \hh_1^*\ |\ \mathfrak{g}_{1, \alpha} \neq (0)\}$. Then $\Delta_{1}^{\times} = \Delta(\gg_1, \hh_1) \backslash \{0\}$ is clearly an irreducible reduced finite
root system having at most two root lengths. Let ${\Delta_{1, \mathrm{sh}}^{\times}}$ denote the set of all non-zero short roots of $\Delta_1$. Define

\[
    {\Delta_{1, \mathrm{en}}^{\times}}= 
\begin{cases}
    \Delta_{1}^{\times} \cup 2  {\Delta_{1, \mathrm{sh}}^{\times}},\,\,\,\,\,\,\,\,\,\,\,\,\,\ &\text{if} \,\,\,\,\Delta_{1}^{\times} \ \text{is \ of \ type} \ B_l \\
    \Delta_{1}^{\times}  \,\,\,\,\,\,\,\,\,\,\,\,\,\,\,\,\,\,\,\,\,\ ,           & \text{otherwise}
\end{cases}
\]\
\,\,\,\,\,\,\,\,\,\,\,\,\,\,\,\,\,\,\,\,\ and $\Delta_{1, \mathrm{en}} = {\Delta_{1, \mathrm{en}}^{\times}} \cup \{0\}$.\\
Suppose that $\{\alpha_1, \cdots ,\alpha_p\}$ is a collection of simple roots of $\mathfrak{g}_1$ with respect to $\mathfrak{h}_1$. 
Then $Q_1 = \bigoplus_{i=1}^{p} \mathbb{Z} \alpha_i$ is the corresponding root lattice containing the non-negative root lattice $Q_1^+ = \bigoplus_{i=1}^{p} \mathbb{Z}_{+} \alpha_i$. Also let  ${\Delta_{1, \mathrm{en}}^{\times,+}}$ and  ${\Delta_{1, \mathrm{en}}^{\times,-}}$ denote the respective set of positive and negative roots of the extended root system.\\
For $\lambda, \mu \in \mathfrak{h}_1^*$, define a partial order relation $``\leqslant_1"$ on $\mathfrak{h}_1^*$ by setting\\
$\lambda \leqslant_1 \mu\ \text{if\ and\ only\ if}\ \mu-\lambda= \sum_{i=1}^{p} n_i \alpha_i\ \text{for\ some}\ n_1, \cdots ,n_p\ \text{in}\ \mathbb{Z}_{+}$.\\
Moreover, if $(\cdot,\cdot)$ is a non-degenerate, symmetric and associative bilinear form on $\mathfrak{g}_1$, then we shall denote the set of all dominant integral weights of $\mathfrak{g}_1$ with respect to $\mathfrak{h}_1$ by $P^{+}_1 = \{ \lambda\ \in \mathfrak{h}_1^*\ | \ (\lambda,\alpha_i) \in \mathbb{Z}_{+}\ \forall \ i=1, \cdots ,p \}$. 

\bdfn\label{Minimal}
A weight $\lambda \in P^{+}_1$ is said to be minimal if for any non-zero weight $\mu \in P^{+}_1$ satisfying $\mu \leqslant_1 \lambda$, we have $\mu = \lambda$.
\edfn

\brmk
Note that zero is always a minimal weight.
\ermk

Let us now define the notion of a Lie torus, which is our main object of interest. In this paper, we shall use the characterization of a Lie torus given in \cite{ABFP} and take it as our definition rather than the axiomatic definition recorded in \cite{Y1,Y2}. For this purpose, we first need the following definition.

\bdfn
A finite-dimensional module $V$ over $\gg_1$ is said to satisfy condition $(M)$ if 
\begin{enumerate}
\item $V$ is an irreducible module having dimension greater than one; 
\item The weights of $V$ relative to $\hh_1$ are contained in $ { \Delta_{1, \mathrm{en}}}$.
\end{enumerate}
\edfn

\bdfn\label{Tori}
A multiloop algebra $L(\mathfrak{g}, \sigma)$ is called a Lie torus
if
\begin{enumerate}
\item $\mathfrak{g}_{\bar{0}}$ is a finite-dimensional simple Lie algebra;
\item For $\bar{k} \neq \overline{0}$ and $\mathfrak{g}_{\bar{k}} \neq (0)$, $\mathfrak{g}_{\bar{k}} \cong U_{\bar{k}} \oplus W_{\bar{k}}$, where $U_{\bar{k}}$
is trivial as an $\gg_{\bar{0}}$-module and either $W_{\bar{k}}$ is zero or it satisfies condition $(M)$;
\item $|\langle\sigma_1, \cdots , \sigma_n \rangle|= \prod_{i=1}^{n} |\sigma_i|$  where $|\sigma_i|$ denotes the order of the automorphism $\sigma_i\ (1 \leqslant i \leqslant n)$ and $|\langle\sigma_1, \cdots , \sigma_n \rangle|$ is the order of the group generated by the $\sigma_i$'s.
\end{enumerate}
\edfn

\brmk\label{Automorphism}
In view of \cite[Theorem 8.5]{K}, it is enough to work with only diagram automorphisms in the one variable setup. Moreover, it was shown in \cite[Proposition 3.4.1]{ABFP} that conditions (1) and (2) (in Definition \ref{Tori}) are  equivalent to the condition that the single automorphism in this case is actually a diagram automorphism. Thus if we now take $\sigma_1$ to be a diagram automorphism and the other automorphisms to be trivial, then it is easy to check that $L(\mathfrak{g}, \sigma)$ also satisfies the above three conditions. The  irreducible integrable representations of these graded multiloop algebras were studied in \cite{BT}.
\ermk
For the rest of this paper, we shall denote the Lie torus $L(\mathfrak{g}, \sigma)$ by $LT$.\\
Let $\mathfrak{h}_{\bar{0}}$ be a Cartan subalgebra of the simple Lie algebra $\mathfrak{g}_{\overline{0}}$. Then it is known that $\mathfrak{h}_{\bar{0}}$ is an ad-diagonalizable subalgebra of $\mathfrak{g}$ \cite[Lemma 3.1.3]{N} and $\Delta^{\times} = \Delta^{\times}(\mathfrak{g}, \mathfrak{h}_{\bar{0}})$ is an irreducible (possibly non-reduced) finite root system \cite[Proposition 3.3.5]{N}. Putting $\Delta_{0} = \Delta(\mathfrak{g}_{\bar{0}},\mathfrak{h}_{\bar{0}})$\ and \ $\Delta = \Delta^{\times} \cup \{0\}$, let us denote the corresponding root latices by $Q_{\overline{0}}$ and $Q$ respectively. Finally let $\Delta_{0}^{+}$ stand for the set of all positive roots of $\mathfrak{g}_{\overline{0}}$ with respect to $\mathfrak{h}_{\overline{0}}$.

\subsection{Properties of LT}\label{Torus} \cite{ABFP}
\begin{itemize}
\item[(LT1)] The Lie torus $LT$ is $\mathbb{Z}^n$-graded as well as $Q$-graded. This is usually referred to as a Lie $\mathbb{Z}^n$-torus of type $\Delta$.
\item[(LT2)]\label{Torus2}
\[\ {\Delta}= 
\begin{cases}
    \Delta_{0, \mathrm{en}}\ , &\text{if} \,\,\,\,\Delta_{0}^{\times} \ \text{is \ of \ type} \ B_l \\
    \Delta_0  \,\,\,\,\,\,\ ,           & \text{otherwise}
\end{cases}
\]
 \item[(LT3)] For each $\overline{k} \in \overline{\Lambda}$, we have 
\begin{align*}
\mathfrak{g}_{\overline{k}} = \bigoplus_{\alpha \in \mathfrak{h}_{\overline{0}}^*} \mathfrak{g}_{\overline{k}}(\alpha)
\end{align*}
where $\mathfrak{g}_{\overline{k}}(\alpha) = \{ x \in \mathfrak{g}_{\overline{k}} \ | \ [h,x] = \alpha(h)x \ \forall\  h \in \mathfrak{h}_{\overline{0}}\}$ and dim$\mathfrak{g}_{\overline{k}}(\alpha) \leqslant1$ for every non-zero $\alpha$ in $\mathfrak{h}_{\overline{0}}^*$.
\item[(LT4)] The centre of $LT$ is trivial and thus $LT$ is referred to as a centreless Lie torus. 
\end{itemize}

\subsection{Universal Central Extension of Lie Torus}\label{Extension}
Let $LT$ be a centreless Lie $\mathbb{Z}^n$-torus of type $\Delta$ as defined earlier. Now just as in the case of $\tau_{A}$ in \eqref{Subspace}, we can similarly construct the spaces $\Omega_{A(\underline{m})}$ and $dA(\underline{m})$ for the smaller algebra $A(\underline{m})= \mathbb{C}[t_1^{\pm m_1}, \cdots, t_n^{\pm m_n}]$. Set $\mathcal{Z}(\underline{m})= \Omega_{A(\underline{m})}/dA(\underline{m})$ and define
\begin{align*}
\overline{LT} := LT \oplus \mathcal{Z}(\underline{m})
\end{align*}
Then $\overline{LT}$ forms a Lie algebra under the following bracket operations.
\begin{enumerate}
\item $[x(\underline{k}), y(\underline{l})] = [x,y](\underline{k} + \underline{l}) \ + \ (x|y) \sum_{i=1}^{n} k_it^{\underline{k} + \underline{l}}K_i$ ;
\item $\mathcal{Z}(\underline{m})$ is central in $\overline{LT}$.
\end{enumerate}
\brmk\label{Generate}
\
\begin{enumerate}
\item Note that since $(\cdot|\cdot)$ is invariant under automorphisms of $\mathfrak{g}$, it is evident that $\underline{k} + \underline{l} \in \Gamma$ whenever $(x|y) \neq 0$ and thus the above bracket operation on $\overline{LT}$ is well-defined. 
\item $\overline{LT}$ is the universal central extension of $LT$ \cite[Corollary 3.27]{J}. We shall sometimes refer to $\mathcal{Z}(\underline{m})$ as the infinite-dimensional centre.
\item $\overline{LT}$ is also a Lie $\mathbb{Z}^n$-torus of type $\Delta$ (with respect to the axiomatic definition of a Lie torus recorded in \cite{Y1}). This implies that $\overline{LT}$ is generated as a Lie algebra by the spaces $(\overline{LT})_{\alpha}$, $\alpha \in \Delta$.
\end{enumerate}
\ermk
As in the case of $LT$, its universal central extension $\overline{LT}$ is also $\mathbb{Z}^n$-graded. To measure this grading, we introduce the derivations $d_1, \cdots ,d_n$ on $\overline{LT}$, the actions of which are given below. Let $D$ denote the linear span of $d_1, \cdots ,d_n$. Set
\begin{equation}\label{Derivations}
\widetilde{LT} = \overline{LT} \oplus D \,\,\,\ \text{and}\,\,\,\ \widehat{LT} = LT \oplus D.
\end{equation}
We now define a Lie bracket on $\widetilde{LT}$ by simply extending
the Lie algebra structure on $\overline{LT}$ to  $\widetilde{LT}$ in the following way.
\begin{enumerate}
\item $d_i(x(\underline{k}))=[d_i, x(\underline{k})] = k_i x(\underline{k})$ ;
\item $d_i(t^{\underline{k}}K_j)=[d_i,t^{\underline{k}}K_j]=k_it^{\underline{k}}K_j$ ;
\item  $[d_i,d_j]=0$. 
\end{enumerate}

\brmk
Note that $Z(\widetilde{LT})= \text{span} \{K_1, \cdots, K_n\}$.
\ermk
\subsection{Roots and Co-roots}\label{Roots}
Consider the subalgebra of $\widetilde{LT}$ given by
\begin{align*}
\widetilde{\mathfrak{h}} = \mathfrak{h}_{\overline{0}} \oplus \displaystyle{\sum_{i =1}^{n}\CC K_i \oplus \sum_{i=1}^{n}} \CC d_i.
\end{align*}
In order to describe the roots of $\widetilde{LT}$, let us first define $\delta_i \in \widetilde{\mathfrak{h}}^*$ by setting
\begin{align*}
\delta_i (\mathfrak{h}_{\overline{0}}) =0,\,\,\,\  \delta_i (K_j) = 0\,\,\,\ \text{and}\,\,\,\ \delta_i (d_j) = \delta_{i j}.
\end{align*}
Put $\delta_{\underline{\beta}} = \displaystyle{\sum_{i =1}^{n}{\beta_i \delta_i}}$ for $\underline{\beta} = (\beta_1, \cdots, \beta_n) \in \mathbb{C}^n$. For $\underline{k} \in \mathbb{Z}^n$, we shall refer to the vector $\delta_{\underline{k} + \underline{\gamma}}$ as the translate of $\delta_{\underline{k}}$ by the vector $\underline{\gamma} = (\gamma_1, \cdots, \gamma_n) \in \mathbb{C}^n$.  Set $\widetilde{\Delta} = \{ \alpha + \delta_{\underline{k}}\ | \ \alpha \in \Delta_{0, \mathrm{en}}^{\times}, \ \underline{k} \in \ZZ^{n} \} \cup \{\delta_{\underline{k}} \ | \ \underline{k} \in \mathbb{Z}^n \setminus \{\underline{0}\} \}$. Then we have the root space decomposition of $\widetilde{LT}$ with respect to $\widetilde{\mathfrak{h}}$ given by
\begin{equation}\label{Decomposition}
\widetilde{LT} = \displaystyle{\bigoplus_{\gamma \in \widetilde{\Delta} \cup \{0\}}}\big ({\widetilde{LT}\big )_{\gamma}}
\end{equation}
where 
\[
 \ {(\widetilde{LT})_{\alpha + \delta_{\underline{k}}}}= 
\begin{cases}
    \gg_{\bar{k}}(\alpha) \otimes \mathbb{C}t^{\underline{k}}\,\,\,\,\,\,\,\,\,\,\,\,\,\,\,\,\,\,\,\,\,\,\,\,\,\,\,\,\,\,\,\,\,\,\,\,\,\,\,\,\,\,\,\,\,\,\,\,\,\,\,\,\,\,\,\,\,\,\,\,\,\,\,\,\,\,\,\,\,\ ,\ \text{if} \,\,\,\,\alpha \neq 0 \\
     \big (\gg_{\bar{k}}(0) \otimes \mathbb{C}t^{\underline{k}}\big ) \oplus  \displaystyle{\big ( \sum_{1 \leq i \leq n, \,\,\underline{k} \in \Gamma \setminus \{\underline{0}\}}}\mathbb{C}{t^{\underline{k}}}K_i \big )\ , \ \text{if}\,\,\,\,\alpha=0,\ \underline{k} \neq \underline{0}\\
      \widetilde{\mathfrak{h}} \,\,\,\,\,\,\,\,\,\,\,\,\,\,\,\,\,\,\,\,\,\,\,\,\,\,\,\,\,\,\,\,\,\,\,\,\,\,\,\,\,\,\,\,\,\,\,\,\,\,\,\,\,\,\,\,\,\,\,\,\,\,\,\,\,\,\,\,\,\,\,\,\,\,\,\,\,\,\,\,\,\,\,\,\,\,\,\,\,\,\,\,\,\,\,\,\,\,\ ,\ \text{if}\,\,\,\,\alpha=0,\ \underline{k}=\underline{0}
\end{cases}
\]
This shows that the roots of $\widetilde{LT}$ are given by $\widetilde{\Delta}$. Furthermore, let
\begin{align*}
\widetilde{\Delta}_{+}= \{ \alpha + \delta_{\underline{k}}\ | \ \alpha \in \Delta_{0, \mathrm{en}}^{\times,+}\ ,\ \underline{k} \in \ZZ^{n} \},\ \ \widetilde{\Delta}_{-}= \{ \alpha + \delta_{\underline{k}}\ | \ \alpha \in \Delta_{0, \mathrm{en}}^{\times,-}\ ,\ \underline{k} \in \ZZ^{n} \}
\end{align*}
denote the set of all positive and negative roots of $\widetilde{LT}$ respectively.\\
A root $\gamma = \alpha + \delta_{\underline{k}}$ is said to be real if $\alpha \neq 0$, otherwise we call it a null root. Let us denote the set of all real roots of $\widetilde{LT}$ by 
$\widetilde{\Delta}^{\mathrm{re}}$. For each such $\gamma \in \widetilde{\Delta}^{\mathrm{re}}$, define the corresponding co-root $\gamma^{\vee} := \alpha^{\vee} + \frac{2}{(\alpha|\alpha)} \displaystyle{\sum_{i =1}^{n}{k_i K_i}}$ where $\alpha^{\vee}$ is the co-root of $\alpha \in \Delta_{0, \mathrm{en}}$. 
\brmk
\
\begin{enumerate}
\item Note that dim$\big ({\widetilde{LT}\big )_{\gamma}} < \infty \ \forall\ \gamma \in \widetilde{\Delta} \cup \{0\}$.
\item $\gamma(\gamma^{\vee}) = \alpha(\alpha^{\vee}) = 2 \ \forall\ \gamma \in \widetilde{\Delta}^{\mathrm{re}}$.
\end{enumerate}
\ermk
\subsection{The Weyl Group}\label{Weyl}
For each $\gamma \in \widetilde{\Delta}^{\mathrm{re}}$, define
the reflection operator $r_{\gamma}$ on $\widetilde{\mathfrak{h}}^*$ by setting
\begin{align*}
 r_{\gamma}(\lambda) = \lambda - \lambda(\gamma^{\vee}) \gamma\ \forall \ \lambda \in \widetilde{\mathfrak{h}}^*.
\end{align*}
Then the Weyl group of $\widetilde{LT}$ which we denote by $W$ is the group
generated by all such reflections $r_{\gamma}$ with $\gamma \in \widetilde{\Delta}^{\mathrm{re}}$. 

\brmk
The Weyl group for the toroidal Lie algebra $\tau_{A(\underline{m})}$ is defined in a similar manner \cite{MS,E3}.
\ermk
\subsection{Change of Coordinates}\label{Twist}
Let $G=GL(n,\mathbb{Z})$ denote the general linear group of $n \times n$ matrices with integer entries and having determinant $\pm 1$. Then $G$ acts on $\mathbb{Z}^n$ via matrix multiplication. Fix any $B \in G$ and put 
\begin{align*}
({d_1}^{\prime},\cdots,{d_n}^{\prime})^T= (B^T)^{-1}({d_1},\cdots,{d_n})^T.
\end{align*}
Taking $\widehat{L}(\mathfrak{g}) = L(\mathfrak{g}) \oplus D$ (see (\ref{Algebra})), let us define an automorphism of $\widehat{L}(\mathfrak{g})$ which we shall denote by $T_{B}$ by setting 
\begin{align*}
T_{B}(x(\underline{k})) = x \otimes t^{B.\underline{k}}\,\,\,\ \text{and} \,\,\,\ T_{B}(d_i)= {d_i}^{\prime}.
\end{align*}
It is straightforward to verify that $T_{B}$ is indeed an automorphism.\\
Set $s_i=t^{B.e_i}$ where $\{e_i\}_{i=1}^n$
denotes the standard basis of $\mathbb{Z}^n$. This shows that $A=  \CC[s_1^{\pm 1}, \cdots, s_n^{\pm 1}]$. Using the  above automorphism $T_{B}$,  we can now define a new graded Lie torus $\widehat{LT}(B)$ by changing the variables from $t_i$ to $s_i$ and the derivations from $d_i$ to ${d_i}^{\prime}$ in $\widehat{LT}$ (see (\ref{Derivations})). Observe that this automorphism takes our original graded Lie torus $\widehat{LT}$ to the newly formed graded Lie torus $\widehat{LT}(B)$.  This phenomenon is called change of coordinates. In this paper, we shall use this change of coordinates without any further comments and simply refer to it as $up$ $to$ $a$ $change$ $of$ $coordinates$.

\brmk\label{Subalgebra} 
Note that by definition, the universal central extension of $\mathfrak{g}_{\bar{0}} \otimes A(\underline{m})$ is simply $\mathfrak{g}_{\bar{0}} \otimes A(\underline{m}) \oplus \mathcal{Z}(\underline{m})$ (see Subsection \ref{Extension}). It is now immediate that the (untwisted) toroidal Lie algebra given by
\begin{align*}
\tau_{A(\underline{m})} = \mathfrak{g}_{\bar{0}} \otimes A(\underline{m}) \oplus \mathcal{Z}(\underline{m}) \oplus D
\end{align*}
is a subalgebra of $\widetilde{LT}$. Strictly speaking, the  space of derivations in this case should also transform accordingly under the automorphism $T_{B}$ (see Subsection \ref{Twist}). But for the sake of notational convenience, we shall stick to our original space of derivations $D$ for this subalgebra $\tau_{A(\underline{m})}$ without any mention of change of coordinates. 
\ermk

\section{Integrable Modules}\label{Representations}
In this section, we define the notion of integrable modules and finally show that we can always find a non-zero highest weight vector in our module.
\bdfn
An $\widetilde{LT}$-module $V$ is called integrable if 
\begin{enumerate}
\item $\displaystyle{V = \bigoplus_{\lambda \in \widetilde{\mathfrak{h}}^*} {V_{\lambda}}}$, where $V_{\lambda} = \{v \in V |\,\, h.v = \lambda(h)v \,\, \forall \,\, h \in \widetilde{\mathfrak{h}}\}$.
\item $\mathrm{dim}(V_{\lambda}) < \infty$ $\forall \ \lambda \in \widetilde{\mathfrak{h}}^*$.
\item All the real root vectors of $\widetilde{LT}$ act locally nilpotently on $V$, i.e. for each $x \in \mathfrak{g}_{\overline{k}}(\alpha)\otimes \mathbb{C}t^{\underline{k}}$ $(\alpha \neq 0)$ and every $v \in V$, there exists $m=m(\alpha,\underline{k},v) \in \mathbb{N}$ such that $x^m.v=0$.
\end{enumerate} 
\edfn
For an integrable module $V$ over $\widetilde{LT}$, let us denote the set of all weights of $V$ by $P(V) = \{\mu \in \widetilde{\mathfrak{h}}^* \ | \ V_{\mu} \neq (0)\}$ . For any $\mu \in P(V)$, $V_{\mu}$ is called the weight space of $V$ having weight $\mu$ and the elements of $V_{\mu}$ are referred to as the weight vectors of $V$ with weight $\mu$.\\
In this paper, our goal is to classify all those irreducible integrable modules over $\widetilde{LT}$ where the elements $K_1, \cdots, K_n$ act trivially on the  module. To serve this purpose, we first need to recall some basic properties of integrable modules that we enlist below.
\bppsn\label{Integrable}
Let $V$ be an integrable module over $\widetilde{LT}$. Then
\begin{enumerate}
\item $P(V)$ is invariant under the action of $W$.
\item $\mathrm{dim}(V_{\lambda}) = \mathrm{dim}(V_{w \lambda})$ $\forall\ \lambda \in P(V)$ and $w \in W$.
\item If $\lambda \in P(V)$  and $\gamma \in \widetilde{\Delta}^{\mathrm{re}}$, then $\lambda (\gamma^{\vee}) \in \ZZ$.
\item If $\lambda \in P(V)$  and $\gamma \in \widetilde{\Delta}^{\mathrm{re}}$ with $\lambda(\gamma^{\vee}) >0$, then $\lambda - \gamma \in P(V)$.
\end{enumerate}
\eppsn

\begin{proof}
The proof is similar to that of \cite[Lemma 2.3]{E3}. The irreducibility of $V$ is not required.
\end{proof}
The root space decomposition of $\widetilde{LT}$ in \eqref{Decomposition} automatically induces a natural triangular decomposition of $\widetilde{LT}$ given by
 \begin{equation}\label{LT}
 \widetilde{LT} = \big (\widetilde{LT}\big )_{-} \oplus \big (\widetilde{LT}\big )_{0} \oplus \big (\widetilde{LT}\big )_{+}
 \end{equation}
 where
\begin{align*}
& \big (\widetilde{LT}\big )_{-} =  \displaystyle{ \bigoplus_{\alpha \in \widetilde{\Delta}_{-}, \,\,\underline{k} \in \ZZ^{n}}}\bigg ({\gg_{\bar{k}}(\alpha) \otimes \mathbb{C}t^{\underline{k}}}\bigg ) ,\\
& \big (\widetilde{LT}\big )_{0} = \bigg (\displaystyle{ \sum_{\underline{k} \in \ZZ^{n}}} {\gg_{\bar{k}}(0) \otimes \mathbb{C}t^{\underline{k}}} \bigg )\ \bigoplus\ \bigg (\displaystyle{ \sum_{1 \leq i \leq n, \,\,\underline{k} \in \Gamma}}\mathbb{C}{t^{\underline{k}}}K_i \bigg ) \ \bigoplus  D,\\
& (\widetilde{LT})_{+} =  \displaystyle{ \bigoplus_{\alpha \in \widetilde{\Delta}_{+}, \,\,\underline{k} \in \ZZ^{n}}}\bigg ({\gg_{\bar{k}}(\alpha) \otimes \mathbb{C}t^{\underline{k}}}\bigg ).
\end{align*}
\brmk\label{Non-abelian}
We would like to emphasize that unlike in the case of toroidal Lie algebras (see (\ref{Abelian})), the portion of the Lie torus $\overline{LT}$ (see Subsection \ref{Extension}) sitting inside $\big (\widetilde{LT}\big )_0$ in our triangular decomposition of $\widetilde{LT}$, namely 
\begin{align*}
\big (\overline{LT}\big )_0 := \bigg (\displaystyle{ \sum_{\underline{k} \in \ZZ^{n}}} {\gg_{\bar{k}}(0) \otimes \mathbb{C}t^{\underline{k}}} \bigg )\ \bigoplus\ \bigg (\displaystyle{ \sum_{1 \leq i \leq n, \,\,\underline{k} \in \Gamma}}\mathbb{C}{t^{\underline{k}}}K_i \bigg )
\end{align*}
modulo $\mathcal{Z}(\underline{m})$ may not be abelian. In fact, it need not even be solvable.
This can be readily checked using \cite[Example 4.3.1]{ABFP} as in this case
\begin{enumerate}
\item $[\gg_{\overline{\lambda}}(0),\gg_{\overline{\mu}}(0)] \neq (0)$\
\text{where} 
$\overline{{\lambda}}=(\overline{1}, \overline{0}, \overline{0}),\
\overline{{\mu}}=(\overline{0}, \overline{1}, \overline{0})$ and 
\item the commutator subalgebra of  $\displaystyle\bigg ({ \bigoplus_{\overline{k} \in \overline{\Lambda} }}\ {\gg_{\bar{k}}(0)\bigg )}$ is a perfect Lie algebra.
\end{enumerate}
\ermk
\blmma\label{L3.4}
\cite[Theorem 3.5]{T} Let $R$ be a ring with unity and $N$ be an $R$-module. Then the  following statements are equivalent:
\begin{enumerate}
\item N is the sum of a family of irreducible submodules.
\item N is the direct sum of a family of irreducible submodules.
\item Every submodule of N is a direct summand.
\end{enumerate}

\elmma
The next lemma is a standard result in the representation theory of finite- dimensional semisimple Lie algebras. Nonetheless, we provide a proof of this statement for the convenience of the reader and also because we are unable to find a precise reference.
\blmma\label{Semisimple}
Let $L$ be a finite-dimensional semisimple Lie algebra with a Cartan subalgebra $H$. Then any integrable $L$-module with respect to $H$ is finite-dimensional.
\elmma

\begin{proof}
As $V$ is an integrable $L$-module, it is immediate that every cyclic submodule of $V$ generated by a weight vector is finite-dimensional. Let us now write 
\begin{align*}
V = \sum_{j \in J}{U(L)v_j}
\end{align*}
where $\{v_j\}_{j \in J}$ is the set of all weight vectors of $V$. Then by Weyl's theorem \cite[Theorem 6.3]{H} and Lemma \ref{L3.4}, it follows that $V$ is a completely reducible $L$-module. This consequently implies that 
\begin{align*}
\displaystyle{V = \bigoplus_{i \in I } {V({\lambda_i})}}
\end{align*}
where $V({\lambda_i})$ is the finite-dimensional irreducible $L$-module with highest weight $\lambda_i$. Thus we only need to show that $I$ is  a finite indexing set. Let us assume the contrary. Now by elementary representation theory of finite-dimensional semisimple Lie algebras, we know that $V(\lambda_{i})$ has a minimal weight for each $i \in I$ (see Definition \ref{Minimal}). For example, in case of $L= \mathfrak{sl}_2(\mathbb{C})$, we know from standard $\mathfrak{sl}_2(\mathbb{C})$ theory that each $V(\lambda_i)$ has either $0$ or $1$ as one of its weights which are precisely the minimal weights of $\mathfrak{sl}_2(\mathbb{C})$. Moreover from Exercise 13 of \cite[Section 13.4]{H}, we also find that, in general, there are only finitely many minimal weights of $L$. Then by the Pigeonhole Principle, at least one of the weight spaces of $V$ has to be infinite-dimensional. This contradiction proves the lemma.
\end{proof}

\blmma\label{L3.6}
Let $V$ be an integrable $\widetilde{LT}$-module. Then there exists some $\lambda \in P(V)$ such that $\lambda + \eta \notin P(V) \ \forall \ \eta \in Q_{\overline{0}}^+ \setminus \{0\}$ where $Q_{\overline{0}}^+$ denotes the non-negative root lattice of $\mathfrak{g}_{\overline{0}}$. 
\elmma

\begin{proof}
By Remark \ref{Subalgebra}, it is clear that $V$ is an integrable $\tau_{A(\underline{m})}$-module with respect to $\widetilde{\mathfrak{h}}$. Now since $Q_{\overline{0}}=Q$, the required result can be obtained by using our Lemma \ref{Semisimple} and by following essentially the same argument presented in \cite[Lemma 2.6]{E3}. Observe that the irreducibility condition in \cite[Lemma 2.6]{E3} is redundant.
\end{proof}
\bppsn\label{Highest}
Let $V$ be an irreducible integrable module over $\widetilde{LT}$ with $K_1, \cdots ,K_n$ acting trivially on $V$. Then there exists a non-zero $v \in V$ such that $\big (\widetilde{LT}\big)_{+}.v=0$.
\eppsn
\begin{proof} In order to prove this proposition, we closely follow the argument given in \cite[Theorem 2.4(ii)]{C} for the affine case. Pick some $\lambda \in P(V)$ as in Lemma \ref{L3.6}. From Proposition \ref{Integrable}, it is clear that $\lambda|_{\mathfrak{h}_{\overline{0}}}$ is a dominant integral weight of $\mathfrak{g}_{\overline{0}}$. We first show that there exists $\mu \in P(V)$ such that $V_{\mu + \alpha + \delta_{\underline{k}}}=(0) \ \forall \ \alpha \in  \Delta_0^{+}, \ \underline{k} \in  \mathbb{Z}^n$. If we have $V_{\lambda + \alpha + \delta_{\underline{k}}}=(0) \ \forall \ \alpha \in  \Delta_0^{+}$ and $\underline{k} \in  \mathbb{Z}^n$, then we are done. So let us assume that there exists $\alpha \in  \Delta_0^{+}$ and $\underline{k} \in  \mathbb{Z}^n$ satisfying $V_{\lambda + \alpha + \delta_{\underline{k}}} \neq (0)$. Putting $\mu = \lambda + \alpha + \delta_{\underline{k}}$, we claim that $V_{\mu + \beta + \delta_{\underline{s}}}=(0) \ \forall \ \beta \in  \Delta_0^{+}, \ \underline{s} \in  \mathbb{Z}^n$. If not, then we have $ V_{\mu + \beta + \delta_{\underline{s}}} \neq (0)$ for some $\beta \in  \Delta_0^{+}$ and $\underline{s} \in  \mathbb{Z}^n$. As $\alpha, \beta \in \Delta_{0}^{+}$, we must have either $(\alpha + \beta | \alpha) > 0$ or $(\alpha + \beta | \beta) > 0$. Without loss of generality, let us suppose that $(\alpha + \beta | \alpha) > 0$ which thereby gives us $(\mu + \beta + \delta_{\underline{s}}| \alpha + \delta_{\underline{k}} + \delta_{\underline{s}}) > 0$. Consequently by Proposition \ref{Integrable}, it immediately follows that $V_{\lambda + \beta } \neq (0)$ which contradicts our initial choice of $\lambda$. Hence the claim. This implies that                
\begin{align*}
V_{\mu + \alpha + \delta_{\underline{k}}}=(0) \ \forall \ \alpha \in  \Delta_0^{+}, \ \underline{k} \in  \mathbb{Z}^n.
\end{align*}
 Again by applying Proposition \ref{Integrable}, we  also get
\begin{align*} 
V_{\mu + 2\alpha + \delta_{\underline{k}}}=(0) \ \forall \ \alpha \in  \Delta_{0}^{+},\ \underline{k}  \in  \mathbb{Z}^n. 
\end{align*}
The proposition now directly follows from (LT2) of Subsection \ref{Torus2}.
\end{proof}

\section{Action of the infinite-dimensional centre}\label{Centre}
For the rest of this paper, unless otherwise stated, $V$ will always stand for an irreducible integrable module for $\widetilde{LT}$ where the elements $K_1, \cdots, K_n$ act trivially. These modules are sometimes referred to as level zero modules. In this section, we shall analyze the action of $\mathcal{Z}(\underline{m})$ on $V$ and finally show that it in fact acts trivially on the entire module.\\
In order to serve this purpose, let us first set
\begin{align*}
V_{+} =\{v \in V \ | \ \big (\widetilde{LT} \big )_{+}.v=0\}
\end{align*}
which is clearly non-zero by Proposition \ref{Highest}. We shall refer to this space as the $highest$ $weight$ $space$ $of$ $V$. The following lemma easily follows from the irreducibility of $V$ and by the PBW Theorem.

\blmma\label{L4.1}\
\begin{enumerate}
\item $V_{+}$ is an irreducible module over $(\widetilde{LT})_{0}$.
\item $V= U \big((\widetilde{LT})_{-}\big)V_{+}$.
\end{enumerate}
\elmma

Note that by our initial assumptions on $V$, it clearly follows that $V_{+}$ is a weight module having finite-dimensional weight spaces with respect to $\widetilde{\mathfrak{h}}$.

\bppsn \label{P4.2}
\cite[Proposition 4.13]{E3} Let $V$ be an irreducible integrable module over $\tau_{A(\underline{m})}$ with respect to $\widetilde{\mathfrak{h}}$ and suppose that $K_1, \cdots ,K_n$ act trivially on $V$. Then $\mathcal{Z}(\underline{m})$ acts trivially on $V$.
\eppsn

\brmk
We cannot use Proposition \ref{P4.2} directly in our case to conclude that $\mathcal{Z}(\underline{m})$ acts trivially on $V$. Because although $V$ is an irreducible module for $\widetilde{LT}$, it may not be irreducible as a $\tau_{A(\underline{m})}$-module.
\ermk

So in order to proceed further, let us consider for each $j = 1, \cdots ,n$,\ the subalgebras of $\tau_{A(\underline{m})}$ given by 
\begin{align*}
L_j = \gg_{\bar{0}} \otimes \mathbb{C}[{t_j}^{m_j},{t_j}^{-m_j}] \bigoplus \mathbb{C}K_j \bigoplus \mathbb{C}d_j.
\end{align*}
Let $W_j$ and $W_{\overline{0}}$ be the Weyl groups related to the affine Kac-Moody algebra $L_j$ and the finite-dimensional simple Lie algebra $\gg_{\bar{0}}$ respectively. Let $\theta$ denote the highest root of $\gg_{\bar{0}}$ and  ${\theta}^{\vee}$ be the  corresponding co-root. Then due to \cite[Proposition 6.5]{K}, we have
\begin{align*}
 W_j \cong W_{\overline{0}} \ltimes T_j
\end{align*}
where $T_j =\{t_{\alpha_{j}} \ | \ \alpha_j \in M_j\}$ and $M_j =\gamma_j(\mathbb{Z} [W_{\overline{0}} \ {\theta}^{\vee}])$
with the natural isomorphism $\gamma_j : \mathfrak{h}_{\overline{0}} \longrightarrow {{\mathfrak{h}}^*_{\overline{0}}}$.

\brmk\label{R4.4}
If we set $\widehat{\mathfrak{h}}_j = \mathfrak{h}_{\overline{0}} \bigoplus \mathbb{C}d_j$, then as $K_j$ acts trivially on $V$, it follows from \cite[Equation 6.5.5] {K} that
\begin{align*}
t_{\alpha_j}(\mu_j) = \mu_j - \mu_j ({\alpha_j}^{\vee}) \delta_j  \ \forall \ \mu_j \in {\widehat{\mathfrak{h}}^*_j},\ \alpha_j \in M_j. 
\end{align*}
\ermk

\blmma\label{L4.5}\
\begin{enumerate}
\item The weights of $V_{+}$ are the same up to a translate of a subset of the null roots, i.e. there exists a unique $\overline{\lambda} \in {{\mathfrak{h}}^*_{\overline{0}}}$ and some $\underline{\beta} \in \mathbb{C}^n$ (not necessarily unique) such that $P(V_{+}) \subseteq \{ \overline{\lambda} + \delta_{\underline{r} + \underline{\beta}} \ | \ \underline{r}  \in \mathbb{Z}^n\}$.
\item $\overline{\lambda}$ is dominant integral.
\end{enumerate}

\elmma

\begin{proof}
(1) Since $\mathfrak{h}_{\overline{0}}$ commutes with  $\big (\widetilde{LT}\big)_{0}$, it follows that $\mathfrak{h}_{\overline{0}}$ acts by scalars on $V_{+}$ by Lemma \ref{L4.1}. As a result, $\mathfrak{h}_{\overline{0}}$ acts by a single linear functional on $V_{+}$ which we shall denote by $\overline{\lambda}$. Moreover the action of the derivations on any non-zero weight vector of $V_+$ is unique upto (translations by) integers.
(2) This follows directly from Proposition \ref{Integrable}.
\end{proof}

\brmk\label{Zero}
\
\begin{enumerate}
\item If $\overline{\lambda}=0$ in Lemma \ref{L4.5}, then using Proposition \ref{Integrable}, we can deduce that the only possible weights of $V$ are of the form $\delta_{\underline{\gamma}}$ $(\underline{\gamma} \in \mathbb{C}^n)$. Then from Remark \ref{Generate}, it follows that $\overline{LT}$ must act trivially on the entire module. Consequently $V$ is isomorphic to the one dimensional highest weight module $V(\delta_{\underline{\gamma}})$ with highest weight $\delta_{\underline{\gamma}}$ for some $\underline{\gamma} \in \mathbb{C}^n$. So without loss of generality, we can assume that $\overline{\lambda} \neq 0$.
\item To simplify our notation, we shall also assume (without any loss of generality) that the vector $\underline{\beta}$ given in Lemma \ref{L4.5} is the null vector.
\end{enumerate}

\ermk

The proof of the following proposition runs almost parallel to the one sketched in \cite[Proposition 3.7]{BES}. However, in this case, we give a much more detailed exposition as it plays a pivotal role in this paper.
\bppsn\label{Finite}
There are only finitely many $\tau_{A(\underline{m})}^0$-submodules of $V_{+}$ where $\tau_{A(\underline{m})} = \tau_{A(\underline{m})}^{-}  \bigoplus \tau_{A(\underline{m})}^0 \bigoplus \tau_{A(\underline{m})}^{+}$ is the triangular decomposition of  the toroidal Lie algebra $\tau_{A(\underline{m})}$ as given in 
(\ref{Abelian}).
\eppsn

\begin{proof}
For each $j = 1, \cdots n$, we have 
$\dfrac{2\theta}{(\theta|\theta)} = \gamma_j(\theta^{\vee}) \in M_j$.\\
Let $\theta_j = \gamma_j(\theta^{\vee})$ and $p_j = (\overline{\lambda}|{\theta_j})\in \mathbb{N}$ (by Lemma \ref{L4.5} and Remark \ref{Zero}).\\
For any $\underline{s} \in \mathbb{Z}^n$, define 
\begin{align*}
\underline{s}(\underline{m}) = (s_1m_1, \cdots , s_nm_n) \in \Gamma.
\end{align*}
Again for any integer $k_j \in \mathbb{Z}$  with $|k_j| \geq p_j$ , there exists non-zero $q_j \in \mathbb{Z}$ and $\ r_j \in \mathbb{Z}_{+}$ such that
\begin{align*}
k_jm_j = q_jp_jm_j + r_jm_j,\ |r_jm_j| \ \textless \ |p_jm_j|
\end{align*}
Set $t_j = t_{\theta_{j}}$ and let $W^{\prime}$ denote the Weyl group corresponding to $\tau_{A(\underline{m})}$.\\ 
Thus if we take
\begin{align*}
w = (\displaystyle \prod_{q_j \textless 0} t_{j}^{-m_jq_j}) (\prod_{q_j \textgreater 0} t_{j}^{m_j q_j}) \in W^{\prime},
\end{align*}
then it is clear from Remark \ref{R4.4} that
\begin{align*}
w(\overline{\lambda} + \delta_{\underline{k}(\underline{m})}) = \overline{\lambda} + \delta_{\underline{r}(\underline{m})},\ |r_im_i| \ \textless \ |p_im_i| \ \forall \ i = 1, \cdots ,n.
\end{align*}
Now set
\begin{align*}
P = \{ \overline{\lambda} + \delta_{\underline{r}(\underline{m})} \ : \underline{r}(\underline{m}) \in \Gamma,   \ |r_im_i| \ \textless \ |p_im_i| \}.
\end{align*}
Finally note that any $\tau_{A(\underline{m})}^0$-submodule of $V_{+}$ is generated by the highest weight vectors of $V$. Hence by virtue of the above argument along with the integrability of $V$ and Remark \ref{Zero}, it follows that any such submodule is in fact generated by a subset of the basis weight vectors of $V$ where all the weights are coming from $P$. But as $P$ is a finite set and the weight spaces of $V$ are given to be finite-dimensional, we get the desired result. 
\end{proof}

\brmk
Note that the proof of the above proposition also shows that the set of dimensions of the weight spaces of $V_{+}$ is uniformly bounded.
\ermk

\bthm\label{Trivial}
$\mathcal{Z}(\underline{m})$ acts trivially on $V$.
\ethm

\begin{proof}
By Proposition \ref{Finite}, we know that there exists a non-zero minimal $\tau_{A(\underline{m})}^0$-submodule of $V_{+}$, say $V_{min}$. Then $V_{min}$ is clearly irreducible. Consider
\begin{align*}
M^{\prime}=U(\tau_{A(\underline{m})})V_{min}
\end{align*}
Since $V_{min}$ is an irreducible $\tau_{A(\underline{m})}^0$-module, it is clear the sum of all the proper $\tau_{A(\underline{m})}$-submodules of $M^{\prime}$ intersects $V_{min}$ trivially. Consequently $V_{min}$ goes injectively to the quotient and hence we get an irreducible quotient of $M^{\prime}$ with finite dimensional weight spaces containing $V_{min}$. Let us denote this $irreducible$ integrable $\tau_{A(\underline{m})}$-module  by $\mathcal{M}$.\\
We can now directly appeal to Proposition \ref{P4.2} to conclude that $\mathcal{Z}(\underline{m})$ acts trivially on $V_{min}$.
Again for any $\underline{r} \in \Gamma$ and $1 \leqslant i \leqslant n,$ it is easy to see that
\begin{align*}
S(\underline{r}, i) = \{ v \in V \ | \ t^{\underline{r}} K_i.\ v = 0\}
\end{align*}
is an $\widetilde{LT}$-submodule of $V$. Thus we are done by the irreducibility of $V$.
\end{proof}

\brmk
The above theorem thereby reduces our final problem to classifying all those irreducible integrable representations of $\widehat{LT}$ (see \eqref{Derivations}) with respect to $\widehat{\mathfrak{h}}$  where $\widehat{\mathfrak{h}} = \mathfrak{h}_{\overline{0}} \oplus D$.
\ermk

\section{Passing on to a module for the multiloop algebra}\label{Non-graded}

In this section, we shall introduce the notion of $highest$ $central$ $operators$ on the highest weight space $V_{+}$. We subsequently use these operators to show that our module $V$ admits of an irreducible finite-dimensional quotient for the multiloop algebra $LT$.\\
Put $\big (LT\big )_0 =  \displaystyle { \bigoplus_{\underline{k} \in \ZZ^{n}}}\bigg({\gg_{\bar{k}}(0) \otimes \mathbb{C}t^{\underline{k}}} \bigg)$ which is a subalgebra of $\widehat{LT}$.
\bdfn\label{Central}
A linear map $z : V_{+} \longrightarrow V_{+}$ is called a highest central operator of degree $\underline{r}$ if it satisfies the following conditions.
\begin{enumerate}
\item $z$ commutes with the action of $\big (LT \big)_0$,  
\item $d_i z -z d_i = r_i z \ \forall \ i = 1, \cdots ,n$.
\end{enumerate}
\edfn

\brmk
\
\begin{enumerate}
\item $\mathfrak{h}_{\overline{0}} \otimes A(\underline{m})$ is an abelian Lie algebra that is central in $\big (LT \big)_0$.
\item For each $h \in \mathfrak{h}_{\overline{0}}$ and $\underline{r} \in \Gamma, \ h \otimes t^{\underline{r}}$ is a highest central operator of degree ${\underline{r}}$ on $V_{+}$.
\end{enumerate} 
\ermk

\blmma\label{Properties}
\
\begin{enumerate}
\item Let $z$ be a highest central operator of degree $\underline{r}$ with $zv \neq 0$ for some $v \in V_{+}$. Then $zw \neq 0$ for every non-zero
$w \in V_{+}$.
\item For a non-zero highest central operator $z$ of degree $\underline{r}$ , there exists a highest central operator $T$ of degree $-\underline{r}$ satisfying $Tz=zT=Id$.
\item The non-zero highest central operators on $V_+$ having the same degree are unique up to scalars.
\end{enumerate} 
\elmma

\begin{proof}
Using our Lemma \ref{L4.1}, the proof follows exactly in the same way as in \cite[Lemma 1.7]{E4} and \cite[Lemma 1.8]{E4}.
\end{proof}

\bppsn\label{Rank}
Let $L = \{ \underline{r} \in \Gamma \ | \ h \otimes t^{\underline{r}} \neq 0 \ on \ V_+ \ for \ some \ h \in \mathfrak{h}_{\overline{0}} \}$ and  $S$ be the subgroup generated by $L$. If we now set $k=rank \ S$, then up to a change of coordinates (see Subsection \ref{Twist}), 
\begin{enumerate}
\item there exist non-zero integers $l_1, \cdots ,l_k$ and non-zero highest central operators $z_1, \cdots ,z_k$ such that the degree of $z_i$ is equal to $l_ie_i$ where $e_i = (0, \cdots,0,1,0, \cdots, 0)\ \forall \ i = 1, \cdots ,k;$
\item $k=n$;
\item $W_+ = \{z_iv\ - v \ | \ v \in V_+,  1 \leqslant i \leqslant k\}$ is a proper $\big (LT \big)_0$-submodule of $V_+$.
\end{enumerate} 
\eppsn
\begin{proof}
(1) First note that $k \neq 0$ by Remark \ref{Zero}. Then the result clearly follows verbatim as in the case of \cite[Theorem 4.5(1)]{E3}. Just observe that the restriction of the automorphism defined in Subsection \ref{Twist} again gives rise to an automorphism of $\mathfrak{h}_{\overline{0}} \otimes A(\underline{m}) \oplus D$.\\
(2) Consider the irreducible $\tau_{A(\underline{m})}$-module $\mathcal{M}$ constructed in Theorem \ref{Trivial}.
Now by \cite[Lemma 3.3]{E2} and \cite[Theorem 3.4]{E2}, there exists a  $\mathbb{Z}^n$-graded algebra homomorphism 
\begin{align*}
\overline{\psi} : U(\mathfrak{h}_{\overline{0}} \otimes A(\underline{m})) \longrightarrow  A(\underline{m})
\end{align*}
and some weight vector $v \in V_{+}$ such that if we put $A_{\overline{\psi}}$ = Im$\overline{\psi}$, then up to a possible change of coordinates, we have
\begin{equation}\label{Full}
U(\mathfrak{h}_{\overline{0}} \otimes A(\underline{m}))v \cong A_{\overline{\psi}} \cong  \CC[t_1^{\pm l_1}, \cdots, t_n^{\pm l_n}]
\end{equation} 
as irreducible $\tau_{A(\underline{m})}^0$-modules. This shows that $S$ indeed has full rank $n$.\\
(3) For each $k =1, \cdots, n-1$, let $D_k$ denote the linear span of $d_{k+1}, \cdots, d_n$. We first show that $W_1 = \{z_1v-v \ | \ v \in V_{+} \}$ is a proper $\big ((LT)_0 \oplus D_1 \big )$-submodule of $V_+$. To this end, consider any non-zero weight vector $v_0$ of $V_+$, say of weight $\mu$. We claim that $v_0 \notin W_1$. If not, then there exist non-zero scalars $c_1, \cdots, c_r$ in $\mathbb{C}$ and non-zero weight vectors $v_1, \cdots, v_r$ in $V_+$ belonging to distinct weight spaces satisfying
\begin{equation}\label{Equation}
v_0 = (c_1 (z_1 v_1) + \cdots c_r (z_1 v_r)) - (c_1 v_1 + \cdots + c_r v_r)
\end{equation}
By means of some elementary arguments involving weights, we can assume (without loss of generality) that 
\begin{align*}
c_i(z_1v_i)=c_{i+1}v_{i+1}\ \forall\ i = 1, \cdots, r-1.
\end{align*}
This eventually reduces our equation \eqref{Equation} to simply
\begin{align*}
v_0 = c_r(z_1v_r)-c_1v_1
\end{align*}
It can be checked that the weights of the above non-zero vectors $z_1v_r$ and $v_1$ have to be distinct. Thus by our initial choice of $v_0$, we have either $c_1=0$ or $c_r=0$,  which is a contradiction. Hence the claim.\\
Let us now consider $W_2 = \{z_2v-v \ | \ v \in V_{+}/W_1 \}$. Again by repeating the above argument, it clearly follows that $W_2$ is a proper $\big ( (LT)_0 \oplus D_2 \big )$-submodule of $V_+/W_1$. Continuing this process finitely many times, we thus get the desired result.
\end{proof}

\brmk\label{R5.5}
Note that from the last proposition, we have $S = \bigoplus_{i=1}^{n}l_i \mathbb{Z}$ up to a change of coordinates.
\ermk
We now prove a key lemma regarding highest central operators which will be crucial in the next section for solving our final classification problem.\\
Set $J(w) = \{X \in U(\mathfrak{h}_{\overline{0}} \otimes A(\underline{m})) \ | \ X.w=0 \}$ for any $w \in V_{+}$.
\blmma\label{Inverse}
Let $z$ be a non-zero highest central operator on $V_+$ of degree $\underline{r} \in S$ and $w$ be any non-zero weight vector in $V_{+}$. Then the unique highest central operator $(modulo \ J(w))$ of degree $-\underline{r}$ on $V_{+}$ given by its inverse lies in $U(\mathfrak{h}_{\overline{0}} \otimes A(\underline{m}))$.
\elmma

\begin{proof} From (\ref{Full}), we have a $\mathbb{Z}^n$-graded isomorphism given by
\begin{align*}
U(\mathfrak{h}_{\overline{0}} \otimes A(\underline{m}))/J(w)  \cong U(\mathfrak{h}_{\overline{0}} \otimes A(\underline{m}))w \cong \CC[t_1^{\pm l_1}, \cdots, t_n^{\pm l_n}]
\end{align*}
This shows that the inverse of any non-zero highest central operator of the form $h \otimes t^{\underline{s}}$ $(h \in \mathfrak{h}_{\overline{0}}$ and $\underline{s} \in \Gamma)$ always lies in $U(\mathfrak{h}_{\overline{0}} \otimes A(\underline{m}))$ modulo $J(w)$. Thus if we now denote the above isomorphism by $\phi$, then we can take 
\begin{align*}
\phi(z.w) = t^{\underline{r}}
\end{align*} 
Consequently we get a unique element  $X(w) \in U(\mathfrak{h}_{\overline{0}} \otimes A(\underline{m}))$ (modulo $J(w)$) such that
\begin{align*}
\phi(X(w).w) = t^{-\underline{r}}
\end{align*}
This element $X(w)$ is clearly of degree $-\underline{r}$. The lemma now follows from the uniqueness of the highest central operator of degree $-\underline{r}$ (upto scalars) on the highest weight space.
\end{proof}

\blmma\label{L5.7}
The number of weights of $V$ are finite up to a translate of a subset of the null roots, i.e. there exist only finitely many $\overline{\lambda_1}, \cdots, \overline{\lambda_{l}} \in \mathfrak{h}_{\overline{0}}^*$ and some $\underline{\beta} \in \mathbb{C}^n$ such that $P(V) \subseteq \{ \overline{\lambda_i} + \delta_{\underline{r} + \underline{\beta}} \ | \ 1 \leqslant i \leqslant l, \ \underline{r} \in \mathbb{Z}^n\}$.
\elmma

\begin{proof}
With regard to Remark \ref{Zero}, we can take 
\begin{align*}
P(V) \subseteq \{ \overline{\lambda} - \eta + \delta_{\underline{r}} \ | \ \eta \in Q_{\overline{0}}^+,\ \underline{r} \in \mathbb{Z}^n \}
\end{align*}
where $\overline{\lambda}$ is the unique element of ${{\mathfrak{h}}^*_{\overline{0}}}$ as in Lemma \ref{L4.5}.\\
Let us now consider any $\overline{\lambda} - \eta + \delta_{\underline{r}} \in P(V)$ for some $\eta \in Q_{\overline{0}}^+$ and $\ \underline{r} \in \mathbb{Z}^n$.\\By Lemma A of \cite[Section 13.2]{H}, there exists $\sigma \in W (\supseteq W_{\overline{0}})$ such that
\begin{align*}
\sigma(\overline{\lambda} - \eta + \delta_{\underline{r}}) =\mu_{\eta} + \delta_{\underline{r}} 
\end{align*}
where $\mu_{\eta}$ is dominant integral. Then from our Proposition \ref{Integrable}, it is clear that $\mu_{\eta} \leqslant \overline{\lambda}$. Moreover from Lemma A of \cite[Section 13.2]{H}, it also follows that $\mu_{\eta} + \delta_{\underline{r}}$ is the unique dominant integral weight conjugate to $\overline{\lambda} - \eta + \delta_{\underline{r}}$ under $W$. But by Lemma B of \cite[Section 13.2]{H}, we know that there can be only finitely many possible $\mu_{\eta}$'s below $\overline{\lambda}$. This proves the lemma.
\end{proof}

\brmk
Combining the above lemma and Proposition \ref{Integrable} along with the Weyl group argument given in Proposition \ref{Finite}, it trivially follows that the set of dimensions of the weight spaces of $V$ is uniformly bounded.
\ermk

\bthm\label{Passing}
Let $V^{\prime} = U(LT)W_+$ and $\overline{V}^{\prime} = V/V^{\prime}$. Then
\begin{enumerate}
\item $V^{\prime}$ is a proper $LT$-submodule of $V$.
\item $\overline{V}^{\prime}$ is a finite-dimensional $LT$-module.
\item $V$ admits of a finite-dimensional irreducible quotient over $LT$, i.e. more precisely, there exists an $LT$-submodule $V^{{\prime}{\prime}}$ of $V$ containing $V^{\prime}$ such that $V/V^{{\prime}{\prime}}$ is a finite-dimensional irreducible module for $LT$.
\end{enumerate}
\ethm

\begin{proof}(1) This is clear from Lemma \ref{L4.1} and Proposition \ref{Rank}.\\
(2) In view of Lemma \ref{L5.7}, it is enough to show that all the weight spaces of $\overline{V}^{\prime}$ are finite-dimensional.
For this purpose, let us first fix some notations.\\
Fix any $\mu \in P(V)$. Put $\overline{\mu} = \mu|_{{\mathfrak{h}}_{\overline{0}}}$ and set
\begin{align*}
\widetilde{V_{\mu}} = \displaystyle {\bigoplus_{\underline{r} \in \ZZ^{n}}} V_{\overline{\mu} + \delta_{\underline{r}}},\,\,\,\,\,\,\,\ \widetilde{V_{\mu}^{{\prime}}} = \widetilde{V_{\mu}} \cap V^{{\prime}}.
\end{align*}
Consequently due to Remark \ref{Zero}, it suffices to show that every weight vector in $\widetilde{V_{\mu}}$ belongs to a finite-dimensional vector space $N(\overline{\mu})$ modulo $\widetilde{V_{\mu}^{{\prime}}}$.\\
From (\ref{Full}), we already know that 
\begin{align}\label{Enveloping}
U(\mathfrak{h}_{\overline{0}} \otimes A(\underline{m}))v  \cong  \CC[t_1^{\pm l_1}, \cdots, t_n^{\pm l_n}]
\end{align}
Choose any non-zero $v_0 \in V_{\overline{\mu} + \delta_{\underline{k}}}$ for some $\underline{k} \in \mathbb{Z}^n$. Then there exist $u_1, \cdots ,u_m \in U(LT)$ and weight vectors $w_1, \cdots ,w_m \in V_{+}$ such that
\begin{equation}\label{Summation}
v_0 = \displaystyle{\sum_{i =1}^{m}} u_i w_i
\end{equation}
Suppose that for each $i=1,\cdots,m$, $w_i$ has weight $\overline{\lambda} + \delta_{\underline{r}^i}$  where $\overline{\lambda}$ is the unique element of ${{\mathfrak{h}}^*_{\overline{0}}}$  obtained in Lemma \ref{L4.5}. Without loss of generality, we can assume that all the above $u_i$'s are weight vectors of $U(LT)$ where $U(LT)$ acts on itself via the action $x.u=xu-ux$.
Then by our choice of $v_0$, the weight of $u_i$ is clearly given by $\overline{\mu}- \overline{\lambda} + \delta_{\underline{k}-\underline{{r}^i}}$.   
Furthermore due to (\ref{Enveloping}), we can choose highest central operators ${z_1}^{\prime}, \cdots, {z_m}^{\prime}$ on $V_{+}$ with respective degrees $\underline{p}^1, \cdots, \underline{p}^m $ (say) such that for each $j=1,\cdots,n$, we have
\begin{align*}
|\ k_j + p^i_{j}\ | < l_j\ \text{where}\ \underline{p}^i = (p^i_{1},\cdots, p^i_{n}).
\end{align*} 
Using \eqref{Summation}, let us now rewrite 
\begin{align*}
v_0 = \displaystyle{\sum_{i =1}^{m}} \bigg ( u_i( w_i - z_i^{\prime}w_i) + u_i (z_i^{\prime}w_i) \bigg )
\end{align*} 
It is now easy to verify that $\displaystyle{\sum_{i =1}^{m}} u_i( w_i - z_i^{\prime}w_i) \in \widetilde{V_{\mu}^{{\prime}}}$. Therefore in order to conclude our proof,  we only need to show that $\sum_{i =1}^{m} u_i (z_i^{\prime} w_i)$ lies in a finite-dimensional subspace of $V$ which is independent of our  initially chosen vector $v_0$. Now observe that if we set $\underline{q}^i = \underline{k} + {\underline{p}^i}$, then
\begin{align*}
u_i (z_i^{\prime}w_i) \in V_{\overline{\mu} + \delta_{\underline{q}^i}}\ \text{where}\ |q^i_{j}| < l_j.
\end{align*} 
But now since the weight spaces of $V$ are finite-dimensional, it is evident that the sum of all the above vectors always lies in a finite-dimensional subspace of $V$, say $N(\overline{\mu})$. Moreover it is also clear that $N(\overline{\mu})$ is independent of the vector $v_0$ and depends only on $\overline{\mu}$. This proves our assertion.\\
(3) For any strictly increasing sequence of $LT$-submodules of $V$ containing $V^{\prime}$ given by
\begin{align*} 
V^{\prime} \subsetneq V_1 \cdots \cdots
\end{align*} 
we clearly obtain a decreasing chain of quotient modules over $LT$ satisfying 
\begin{align*}
0 \neq  dim(V/V_1) \ \textless \ dim(V/V^{\prime})
\end{align*} 
Then by (2), the above chain of increasing submodules must terminate. Thus there exists a proper maximal submodule of $V$ containing $V^{\prime}$, say $V^{{\prime}{\prime}}$. Henceforth we shall denote this irreducible quotient module $V/V^{{\prime}{\prime}}$ by $\overline{V}$.
\end{proof}

\section{Recovering the original module}\label{Graded}
In this section, our aim is to recover the irreducible $\widehat{LT}$-module $V$ from the finite-dimensional irreducible $LT$-module $\overline{V}$ that we have just defined at the end of the previous section. The ideas used in this context are inspired from \cite{E4}, but the arguments applied in that case need to be suitably modified.\\
To this end, first recall that $LT$ is naturally $\mathbb{Z}^n$-graded.
Using this natural gradation, let us now define an $\widehat{LT}$-module structure $\rho(\underline{\alpha})$ on $\overline{V} \otimes A$ where $\underline{\alpha} = (\alpha_1, \cdots, \alpha_n) \in \mathbb{C}^n$ by setting
\begin{align}\label{Loop}
X(\underline{r}).(\overline{v} \otimes t^{\underline{s}}) & =   (X(\underline{r})\overline{v}) \otimes t^{\underline{r} + {\underline{s}}} ,\ X(\underline{r}) \in (LT)_{\underline{r}},\ \overline{v} \in \overline{V},\ \underline{r}, \underline{s} \in \mathbb{Z}^n;\\
d_i.(\overline{v} \otimes t^{\underline{s}})  & = (\alpha_i + s_i) (\overline{v} \otimes t^{\underline{s}}) \ \forall \ 1 \leqslant i \leqslant n.
\end{align}

\brmk\label{Weight}
Observe that any $\lambda \in \mathfrak{h}_{\overline{0}}^*$ can be extended to $\widehat{\mathfrak{h}}^*$ by simply setting $\lambda|_{D}=0$. Then $\overline{V} \otimes A$ is clearly an integrable $\widehat{LT}$-module under the above action with respect to $\widehat{\mathfrak{h}}$ and $( \overline{V} \otimes A )_{\lambda + \delta_{\underline{r}+ \underline{\alpha}}} = \overline{V}_{\lambda} \otimes  \mathbb{C}t^{\underline{r}} \ \forall \ \underline{r} \in \mathbb{Z}^n$.
\ermk

For any $\overline{v} \in \overline{V}$ and $\underline{r} \in \mathbb{Z}^n$, define $\overline{v}(\underline{r}):= \overline{v} \otimes t^{\underline{r}}$.
Next observe that for any fixed choice of $\mu \in P(V_{+})$, the action of the derivations $d_1, \cdots ,d_n$ on the weight $\mu$ directly gives rise to a $\mathbb{Z}^n$-grading on $V$ as an $\widehat{LT}$-module, say
\begin{align*}
V = \displaystyle {\bigoplus_{\underline{r} \in \ZZ^{n}}} V_{{\underline{r}}}
\end{align*} 
The above gradation on $V$ then induces a natural map
\begin{equation} \label{Map}
\phi : V \longrightarrow \overline{V} \otimes A
\end{equation}
$\,\,\,\,\,\,\,\,\,\,\,\,\,\,\,\,\,\,\,\,\,\,\,\,\,\,\,\,\,\,\,\,\,\,\,\,\,\,\,\,\,\,\,\,\,\,\,\,\,\,\,\,\,\,\,\,\,\,\,\,\,\,\,\,\,\,\,\,\,\,\,\,\,\,\,\,\,\,\,\,\,\,\,\,\ v \longmapsto \ \overline{v} ({\underline{r}}),\ v \in V_{{\underline{r}}}$\\
It is trivial to check that $\phi$ is an $\widehat{LT}$-module homomorphism.\\ 
Moreover by Theorem \ref{Passing}, it is also evident that $\phi$ cannot possibly vanish on any non-zero weight vector of $V$. Therefore by the irreducibility of $V$, it follows that $\phi$ must be $injective$.
This proves that $\widehat{V}:= \phi(V)$ is a non-zero $irreducible$ $\widehat{LT}$-submodule of $\overline{V} \otimes A$.\\
Fix any $\mu \in P(V_+)$ and pick a non-zero vector $v$ in $V_{\mu} \cap V_+$.\\
For each $i = 1, \cdots ,n$, set $\alpha_i = \mu(d_i)$ and set  $\underline{\alpha}=
(\alpha_1, \cdots, \alpha_n) \in \mathbb{C}^n$.
With respect to this choice, $v$ is a graded vector in $V$ with $v \in V_{\underline{0}}$. In order to keep track of this grading, let us denote this vector $v$ by $v(\underline{0})$ and its corresponding image $\phi(v)$ by $\overline{v}(\underline{0})$.
Using Theorem \ref{Passing}, it is trivial to check that $\overline{v}(\underline{0}) \neq 0$ which finally shows that $\widehat{V} = U(\widehat{LT})(\overline{v}(\underline{0}))$.\\
Now for any $\underline{s} \in \mathbb{Z}^n$, let us define an $\widehat{LT}$-module structure on $\mathbb{C}$ by setting
\begin{align*}
X(\underline{r}).1 & =  0 \ \forall \ X \in (LT)_{\underline{r}},\ \underline{r} \in \mathbb{Z}^n ;\\
d_i.1  & = -s_i \ \forall \ i = 1, \cdots ,n.
\end{align*}
This clearly gives rise to an $isomorphism$ of $\widehat{LT}$-modules via the map
\begin{align*}
\psi : U(\widehat{LT})(\overline{v}(\underline{0})) \longrightarrow  U(\widehat{LT})(\overline{v}(\underline{s})) \otimes \mathbb{C}\\
\overline{w} ({\underline{r}}) \longmapsto \overline{w} ({\underline{r} + \underline{s}}) \otimes 1, \ \overline{w} \in \overline{V}, \underline{r} \in \mathbb{Z}^n
\end{align*}
which shows that $U(\widehat{LT})(\overline{v}(\underline{s}))$ is an $irreducible$ $\widehat{LT}$-submodule of $\overline{V} \otimes A$.

\brmk
\
\begin{enumerate}
\item $\psi$ is always an isomorphism between $ U(\widehat{LT})(\overline{v}(\underline{0}))$ and $U(\widehat{LT})(\overline{v}(\underline{s}))$ as $LT$-modules, but need not be an isomorphism of $\widehat{LT}$-modules.
\item In the above setting, $U(\widehat{LT})(\overline{v}(\underline{r}))$ and $U(\widehat{LT})(\overline{v}(\underline{s}))$ are said to be isomorphic (as $\widehat{LT}$-modules) up to a grade shift for all $\underline{r}, \underline{s} \in \mathbb{Z}^n$. In this case, the derivation $d_i$ acts by scalar multiples of  $(r_i-s_i)$ on $\mathbb{C}$.
\end{enumerate}
\ermk

\bppsn\label{P6.3}
\
\begin{enumerate}
\item $U(\widehat{LT})(\overline{v}(\underline{r}))$ is an irreducible $\widehat{LT}$-module for all $\underline{r} \in \mathbb{Z}^n$.
\item $\overline{V} \otimes A = \displaystyle{\sum_{\underline{r} \in \mathbb{Z}^n}}U(\widehat{LT})(\overline{v}(\underline{r}))$.
\item $\overline{V} \otimes A = \displaystyle{\sum_{0 \leqslant r_i < l_i}}U(\widehat{LT})(\overline{v}(\underline{r}))$.
\end{enumerate}
\eppsn
\begin{proof} (1) Clear from the previous discussion.\\
(2) Let $\overline{w}(\underline{r}) \in \overline{V} \otimes A$ where 
$\overline{w} \in \overline{V}$ and $\underline{r} \in \mathbb{Z}^n$.
Now since we already know that $\overline{V}$ is an irreducible $LT$-module, there exists some $X \in U(LT)$ such that $X.\overline{v} = \overline{w}$. Again observe that the $\mathbb{Z}^n$-grading on $LT$ directly induces a $\mathbb{Z}^n$-grading on $U(LT)$ and therefore there exist $X_1, \cdots, X_m \in U(LT)_{\underline{p}^i}$ such that $X = \displaystyle{\sum_{i=1}^{m}} X_i$. Then for our particular choice of $\underline{\alpha}$ (corresponding to our chosen vector $v$), the action of $LT$ on $\overline{V} \otimes A$ as defined in \eqref{Loop} shows that
\begin{align*}
 \displaystyle{\sum_{i=1}^{m}} X_i.  \overline{v}(\underline{r} - {\underline{p}^i}) = {\sum_{i=1}^{m}}( X_i.  \overline{v})(\underline{r}) = \overline{w}(\underline{r}).
\end{align*} 
(3) By Remark \ref{R5.5}, it suffices to show that for any fixed $\underline{k} \in \mathbb{Z}^n$,
\begin{align*}
U(\widehat{LT})(\overline{v}(\underline{r})) = U(\widehat{LT})(\overline{v}(\underline{r} + \underline{k})) \ \forall \ \underline{r} \in S.
\end{align*} 
To this end, let us consider any non-zero $\underline{r} \in S$. Then by Remark \ref{R5.5}, $\underline{r} = \sum_{i=1}^{n}b_ie_i$ for some $b_1, \cdots ,b_n \in \mathbb{Z}$. Now set
\begin{align*}
z = \displaystyle \prod_{i=1}^n z_i^{b_i}
\end{align*} 
This immediately implies that $zw-w \in W_+ \ \forall \ w \in V_+$. In particular, we have $zv-v \in W_+$ by our initial choice of $v$.
And so from our Lemma \ref{Inverse}, we can conclude that
\begin{align*}
\overline{v}(\underline{r} + \underline{k}) = (\overline{z.v})(\underline{r} + \underline{k}) = z.(\overline{v}(\underline{r})) \in U(\widehat{LT})(\overline{v}(\underline{r}))
\end{align*} 
which proves the proposition by (1).
\end{proof}
Since the $\widehat{LT}$-module homomorphism $\phi$ defined in (\ref{Map}) is known to be injective, the following theorem is now a direct  consequence of the above proposition and Lemma \ref{L3.4}.
\bthm\label{Reducible}
\
\begin{enumerate}
\item $(\overline{V} \otimes A, \ \rho(\underline{\alpha}))$ is a completely reducible integrable module for $\widehat{LT}$ with respect to $\widehat{\mathfrak{h}}$. The number of irreducible components occurring in this decomposition are finite with every irreducible component being isomorphic to each other up to a grade shift.
\item For a suitable choice of $\underline{\alpha} \in \mathbb{C}^n, \ V$ is isomorphic to an irreducible component of $(\overline{V} \otimes A, \ \rho(\underline{\alpha}))$ occurring in the above decomposition as $\widehat{LT}$-modules.
\end{enumerate}
\ethm

\section{The final classification theorem}\label{Final}

The last two sections clearly reduce our study of irreducible integrable modules for $\widehat{LT}$ to finite-dimensional irreducible modules over $LT$. These  finite-dimensional irreducible modules have been completely classified in \cite{L}. In this section, we shall recall the classification theorem given in \cite{L} and henceforth classify all the irreducible integrable modules for $\widetilde{LT}$ up to isomorphism. Throughout this section, we shall denote the set of all non-zero complex numbers by $\mathbb{C}^{\times}$.\\
Let $a_{ij} \in \mathbb{C} \ (1 \leqslant i \leqslant k, 1 \leqslant j \leqslant n)$ be such that for each $1 \leqslant i \leqslant k$, we have
\begin{align*}
\underline{a_i} = (a_{i1}, \cdots ,a_{in}) \in (\mathbb{C}^{\times})^n.
\end{align*} 
Again for each $1 \leqslant j \leqslant n  $, choose $b_{ij} \in \mathbb{C}$ such that $b_{ij}^{m_j}=a_{ij}$. Finally set
\begin{align*}
\underline{b_i} = (b_{i1}, \cdots ,b_{in}), \ \ \underline{m} (\underline{b_i}) = (b_{i1}^{m_1}, \cdots ,b_{in}^{m_n})\ \forall \ 1 \leqslant i \leqslant k.
\end{align*} 
Let $V(\lambda_i) \ (1 \leqslant i \leqslant k)$ be the finite-dimensional irreducible module over $\mathfrak{g}$ having highest weight $\lambda_i$. We  shall write $V(\lambda_i, \underline{b_i})$ for the $LT$-module $V(\lambda_i)$ equipped with the following action of $LT$.
\begin{align}\label{Evaluation}
(x \otimes f(t)).v = f(\underline{b_i})x.v
\end{align}
In the literature, these modules are referred to as evaluation modules.\\
We shall denote the $k$-fold tensor product of such a family of evaluation modules by
\begin{align*}
V(\underline{\lambda}, \underline{b}, k) = V(\lambda_1, \underline{b_1}) \otimes \cdots \otimes V(\lambda_k, \underline{b_k})
\end{align*}
with $\underline{b} = (\underline{b_1}, \cdots , \underline{b_k}) \in \mathbb{C}^{nk}$ and $\underline{\lambda} = (\lambda_1, \cdots ,\lambda_k)$ where $\lambda_1, \cdots ,\lambda_k$ are dominant integral weights of $\mathfrak{g} $.\\
We can now give an $LT$-module structure on $V(\underline{\lambda}, \underline{b}, k)$ by simply extending the action of $LT$ given in \eqref{Evaluation} to the whole space, i.e. more precisely,
\begin{align*}
(x \otimes f(t)).(v_1 \otimes \cdots \otimes v_k) = \sum_{i=1}^{k} f(\underline{b_i})(v_1 \otimes \cdots \otimes x.v_i \otimes \cdots \otimes v_k).
\end{align*}
It was shown in \cite[Theorem 4.12]{L} that $V(\underline{\lambda}, \underline{b}, k)$ is a finite-dimensional irreducible $LT$-module if $\underline{m} (\underline{b_i}) \neq \underline{m} (\underline{b_j}) \ \forall \ 1 \leqslant i \neq j \leqslant k$.
\bppsn\label{Multiloop}
\cite[Corollary 4.11]{L} Suppose that $V$ is a finite-dimensional irreducible module over the multiloop algebra $LT$. Then there always exist $\underline{b_1}, \cdots , \underline{b_k} \in (\mathbb{C}^{\times})^n$ and dominant integral weights $\lambda_1, \cdots ,\lambda_k$ of $\mathfrak{g}$ for some $k \in \mathbb{N}$ such that $V \cong V(\underline{\lambda}, \underline{b}, k)$ where $\underline{m} (\underline{b_i}) \neq \underline{m} (\underline{b_j}) \ \forall \ 1 \leqslant i \neq j \leqslant k$.
\eppsn
We can now define an $\widehat{LT}$-module structure $\rho(\underline{\alpha})$ on $V(\underline{\lambda}, \underline{b}, k) \otimes A$ where $\underline{\alpha} = (\alpha_1, \cdots, \alpha_n) \in \mathbb{C}^n$ by setting
\begin{align*}
X (\underline{r}).(v \otimes t^{\underline{s}}) & =   (X(\underline{r})v) \otimes t^{\underline{r} + {\underline{s}}},\ X(\underline{r}) \in (LT)_{\underline{r}},\ v \in V(\underline{\lambda}, \underline{b}, k),\ \underline{r}, \underline{s} \in \mathbb{Z}^n;\\
d_i.(v \otimes t^{\underline{s}})  & = (\alpha_i + s_i) (v \otimes t^{\underline{s}}) \ \forall \ i = 1, \cdots ,n.
\end{align*}
In view of Remark \ref{Weight}, it is clear that $(V(\underline{\lambda}, \underline{b}, k) \otimes  A, \rho(\underline{\alpha}))$ is an integrable $\widehat{LT}$-module. By Theorem \ref{Reducible}, it also follows that $(V(\underline{\lambda}, \underline{b}, k) \otimes  A, \rho(\underline{\alpha}))$ is a completely reducible module over $\widehat{LT}$ with only finitely many irreducible components. Moreover all these irreducible components occurring in this decomposition are mutually isomorphic up to a grade shift. Again due to Proposition \ref{P6.3} and Theorem \ref{Trivial}, we can further conclude that each such irreducible component of $V(\underline{\lambda}, \underline{b}, k) \otimes A$ is in fact isomorphic to an irreducible integrable module over $\widetilde{LT}$ where the elements $K_1, \cdots ,K_n$ act trivially. In this paper, we have shown that these irreducible components completely exhaust all such irreducible integrable modules (see Theorem \ref{Reducible}). Thus we have finally proved the following theorem.

\bthm\label{Classification}
Let $V$ be an irreducible integrable $\widetilde{LT}$-module and suppose that all the central elements $K_1, \cdots ,K_n$ act trivially on $V$. Then up to a possible change of coordinates, $V$ is isomorphic to an irreducible component of $(V(\underline{\lambda}, \underline{b}, k) \otimes A,\ \rho(\underline{\alpha}))$ for some $k \in \mathbb{N}, \ \underline{\lambda} = (\lambda_1, \cdots ,\lambda_k)$, $\underline{\alpha} \in \mathbb{C}^n$ and $\underline{b} = (\underline{b_1}, \cdots , \underline{b_k}) \in \mathbb{C}^{nk}$ where $\underline{b_i} \in (\mathbb{C}^{\times})^n$ for each $1 \leqslant i \leqslant k$ such that $\underline{m} (\underline{b_i}) \neq \underline{m} (\underline{b_j}) \ \forall \ 1 \leqslant i \neq j \leqslant k$, with all the $\lambda_i$'s being dominant integral weights of $\mathfrak{g}$.
\ethm

\brmk 
From Remark \ref{Automorphism} and Theorem \ref{Trivial}, it is completely clear that our final theorem also recovers the classification theorem for irreducible integrable modules having finite-dimensional weight spaces for the graded multiloop algebras considered in \cite{BT}.
\ermk

\section{Afterword}
The following problem appears quite natural and might be addressed in a subsequent paper. \\\\
$\textbf{Problem}.$ Let $V_1$ and $V_2$ be level zero irreducible integrable $\widetilde{LT}$-modules which are irreducible components of $\overline{V}_1 \otimes A$ and $\overline{V}_2 \otimes A$ respectively where $\overline{V}_1$ and $\overline{V}_2$ are finite-dimensional irreducible $LT$-modules. If $V_1 \cong V_2$ as $\widetilde{LT}$-modules, then is $\overline{V}_1 \cong \overline{V}_2$ as $LT$-modules?

\end{document}